\documentstyle[11pt,leqno]{article}
\newtheorem{theorem}{{\sc Theorem}}
\newcommand{\bt}{\begin{theorem}}
\newcommand{\et}{\end{theorem}}
\newcommand{\newsection}[1]{\setcounter{equation}{0} \setcounter{theorem}{0}
\section{#1}}

\newcommand{\NI}{\noindent}
\newcommand{\bea}{\begin{eqnarray}}
\newcommand{\eea}{\end{eqnarray}}

\def \spec#1 {\mathop{#1}}

\def \b #1 {\bf #1}

\newcommand {\CC}{\centerline}

\newcommand{\clf}{{\cal F}}

\newcommand{\cly}{{\cal Y}}
\newcommand{\clz}{{\cal Z}}

\newcommand{\ity}{\infty}
\newcommand{\raro}{\rightarrow}

\newcommand{\vsp}{\vskip 1em}

\newcommand{\be}{\begin{equation}}
\newcommand{\ee}{\end{equation}}
\newcommand{\ben}{\begin{eqnarray*}}
\newcommand{\een}{\end{eqnarray*}}

\begin{document}
\sloppy
\CC {\Large{\bf  Maximum Likelihood Estimation for Stochastic Differential}}
\CC{\Large{\bf Equations Driven by a Mixed Fractional Brownian Motion}} 
\CC{\Large{\bf with Random Effects}}
\vsp
\CC {\bf B.L.S. Prakasa Rao}
\CC{CR RAO Advanced Institute of Mathematics, Statistics}
\CC{ and Computer Science, Hyderabad, India}
\CC{email: blsprao@gmail.com}
\vsp
\NI{\bf Abstract:} We discuss maximum likelihood estimation of parameters for  models governed by a stochastic differential equation driven by a mixed fractional Brownian motion with random effects. 
\vsp 
\NI{\bf Keywords :} Stochastic differential equation; Random effects; Maximum likelihood estimation;  Mixed fractional Brownian motion.  
\vsp 
\NI {\bf AMS Subject Classification:} 60G22, 62M09. 

\newsection{Introduction}

Stochastic modeling by processes driven by a fractional Brownian motion has been used for  phenomena with long range dependence. Statistical inference  for stochastic processes satisfying stochastic differential equations driven by a fractional Brownian motion has been studied earlier and a comprehensive survey of various methods is given in Mishura [1] and in Prakasa Rao [2]. However it was observed that modeling of the financial markets by processes driven by fractional Brownian motion lead to arbitrage opportunities which is contrary to the rational market behaviour. Cheridito [3] proposed modeling by processes driven by a mixed fractional Brownian motion to avoid this problem. There has been a recent interest to study problems of statistical inference for stochastic processes driven by a mixed fractional Brownian motion (mfBm). Existence and uniqueness for solutions of stochastic differential equations driven by a mfBm are investigated in Mishura and Shevchenko [4], Shevchenko [5],  Guerra and Nualart [6] and more recently by Luis da Silva et al. [7] among others . Maximum likelihood estimation for estimation of drift parameter in a linear stochastic differential equations driven by a mfBm is investigated in Prakasa Rao [8]. The method of instrumental variable estimation for such parametric models is investigated in Prakasa Rao [9]. Some applications of such models in finance are presented in Prakasa Rao [10,11]. Stochastic differential equation models with random effects are used in the biomedical field for the study of repeated measurements collected on a series of individuals/subjects. For instance, see Antic et al. [12], Delattre et al. [13], Ditlevsen and  De Gaetano [14], Nie and Yang [15], Nie [16,17], Picchini et al. [18] and Picchini and Ditlevsen [19] among others. 
Parametric inference for linear stochastic differential equations driven by a mixed fractional Brownian motion with random effects based on discrete observations has been studied in Prakasa Rao [31]. Nonparametric estimation for fractional diffusion processes with random effects has been investigated in El Omari et al. [20]. They study the properties of kernel and histogram estimators for estimation of the density of random effects. We discussed nonparametric estimation for models governed by stochastic differential equations with random effects driven by a mixed fractional Brownian motion (mfBm) with Hurst index $H >\frac{1}{2}$  in Prakasa Rao [21]. For parametric inference for processes driven by mfBm, see Marushkevych [22], Rudomino-Dusyatska [23], Song and Liu [24], Mishra and Prakasa Rao [25], Prakasa Rao [26] and Miao [27] among others. El Omari et al. [28] studied estimation of parameters $\mu$ and $\sigma^2$ when the random effects are Gaussian with mean $\mu$ and variance $\sigma^2$ based on discrete observations on the process. Maximum likelihood estimation of stochastic differential equations with random effects driven by fractional brownian motion is investigated in Dai et al. [32]. Our aim in this paper is to  study similar problem for processes driven by mfBm.
\vsp
\newsection{Mixed fractional Brownian motion}
We will now summarize some properties of stochastic processes which are solutions of stochastic differential equations driven by a mixed fractional Brownian motion . We assume that sufficient conditions hold to ensure existence and uniqueness of the solution.
\vsp
Let $(\Omega, \clf, (\clf_t), P) $ be a stochastic basis
satisfying the usual conditions. The natural filtration of a
stochastic process is understood as the $P$-completion of the
filtration generated by this process. Let $\{W_t, t \geq 0\}$ be a standard Wiener process and $W^H= \{W_t^H, t \geq 0 \}$ be an {\it independent}  normalized  fractional Brownian motion with Hurst parameter $H \in (0,1)$, that is, a Gaussian process with continuous sample paths such that $W_0^H=0, E(W_t^H)=0$ and
\be
E(W_s^H W_t^H)= \frac{1}{2}[s^{2H}+t^{2H}-|s-t|^{2H}], t \geq 0, s \geq 0.
\ee
Let
$$\tilde W_t^H= W_t+ W_t^H, t \geq 0.$$
The process $\{\tilde W_t^H, t \geq 0\}$ is called the mixed fractional Brownian motion with Hurst index $H.$ We assume here after that Hurst index $H$ is known. Following the results in Cheridito [3], it is known that the process $\tilde W^H$ is a semimartingale in its own filtration if and only if either $H=1/2$ or $H \in (\frac{3}{4},1].$ 
\vsp
Let us consider a stochastic process $Y=\{Y_t, t \geq 0\}$ defined
by the stochastic integral equation \be Y_t= \int_0^t C(s) ds  +\tilde W_t^H, t \geq 0 \ee where the process $C=\{C(t), t
\geq 0\}$ is an $(\clf_t)$-adapted process. For convenience, we write
the above integral equation in the form of a stochastic
differential equation
\be
dY_t= C(t) dt + d\tilde W_t^H, t \geq 0
\ee
driven by the mixed fractional Brownian motion $\tilde W^H.$ Following the recent works by Cai et al. [29] and Chigansky and Kleptsyna [30], one can construct an integral transformation that transforms the mixed fractional Brownian motion $\tilde W^H$ into a martingale $M^H.$ Let $g_H(s,t)$ be the solution of the integro-differential equation
\be
g_H(s,t)+H \frac{d}{ds}\int_0^t g_H(r,t)|s-r|^{2H-1} sign(s-r)dr=1, 0<s<t.
\ee
Cai et al. [28] proved that the process
\be
M_t^H= \int_0^tg_H(s,t)d\tilde W_s^H, t \geq 0
\ee
is a Gaussian martingale with quadratic variation
\be
<M^H>_t= \int_0^tg_H(s,t)ds, t \geq 0
\ee
For convenience, we denote the function $<M^H>_t$ by $w^H_t.$ Furthermore the natural filtration of the martingale $M^H$ coincides with that of the mixed fractional Brownian motion $\tilde W^H.$ Suppose that, for the martingale $M^H$ defined by the equation (2.5), the sample paths of the process $\{C(t), t \geq 0\}$ are smooth enough in the sense that the process
\be
Q_H(t)= \frac{d}{d<M^H>_t}\int_0^tg_H(s,t)C(s)ds, t \geq 0
\ee
is well defined. Define the process
\be
Z_t= \int_0^tg_H(s,t) dY_s, t \geq 0.
\ee
As a consequence of the results in Cai et al. [29], it follows that the process $Z$ is a fundamental semimartingale associated with the process $Y$ in the following sense.
\vsp
\NI{\bf Theorem 2.1:} {\it Let $g_H(s,t)$ be the solution of the equation (2.4). Define the process $Z$ as given in the equation (2.8). Then the following relations hold.

\NI(i) The process $Z$ is a semimartingale with the decomposition
\be
Z_t= \int_0^tQ_H(t)d<M^H>_s + M^H_t, t \geq 0
\ee
where $M^H$ is the martingale defined by the equation (2.5).

\NI(ii) The process $Y$ admits the representation
\be
Y_t=\int_0^t\hat g_H(s,t)dZ_s, t \geq 0
\ee
where
\be
\hat g_H(s,t)= 1-\frac{d}{d<M^H>_s}\int_0^tg_H(r,s)dr.
\ee
\NI (iii) The natural filtrations $(\cly_t)$ and $(\clz_t)$ of the processes $Y$ and $Z$ respectively coincide.}
\vsp
Applying Corollary 2.9 in Cai et al. [29], it follows that the probability measures $\mu_Y$ and $\mu_{\tilde W^H}$ generated by the processes $Y$ and $\tilde W^H$ on an interval $[0,T]$ are absolutely continuous   with respect to each other and the Radon-Nikodym derivative is given by
\be
\frac{d\mu_Y}{d\mu_{\tilde W^H}}(Y)= \exp[\int_0^TQ_H(s) dZ_s-\frac{1}{2}\int_0^T[Q_H(s)]^2 d<M^H>_s]
\ee
which is also the likelihood function based on the observation $\{Y_s,0\leq s \leq T.\}$ Since the filtrations generated by the processes $Y$ and $Z$ are the same, the information contained in the families of $\sigma$-algebras  $(\cly_t)$ and $(\clz_t)$  is the same and hence the problem of the estimation of the parameters involved  based on the observation $\{Y_s, 0\leq s \leq T\}$ and $\{Z_s, 0\leq s \leq T\}$ are equivalent.
\vsp
\newsection{Preliminaries}

Let us consider a system of stochastic differential equations
\be
dX^i_t= S(X^i_t,\phi_i )\; dt + d \tilde W^{H,i}_t, X^i_0=x^i , 0 \leq t \leq 1, i=1,\dots, N 
\ee
where the processes $\tilde W^{H,i}, 1 \leq i \leq N $ are {\it independent}  mixed fractional Brownian motions with  common Hurst index $H > \frac{1}{2}.$ Suppose the random effects $\phi_i, i=1,\dots,N$ are independent and identically distributed as that of the distribution of $\phi.$ We assume that the random variables $\phi_1,\dots, \phi_N$ are {\it independent} of the random processes $\{\tilde W^{H,i}, i=1,\dots,N\}$ and sufficient conditions hold for the existence and uniqueness of the solutions for the system (3.1). See Mishura and Shevchenko [4], Shevchenko [5], Guerra and Nualart [6] and Luis da Silva [7] for sufficient conditions on the function $S(.,.)$ for the existence and uniqueness of the solution  for the system (3.1).    The problem is to estimate the unknown parameters based on the set of observations $\{ X^i(t), 0\leq t \leq T; i=1,\dots,N \}$ and study the asymptotic properties of the estimators as $N$ tends to infinity. 
\vsp
Here after we assume that the random effect $\phi $ has a probability density function  $g(\psi,\theta)$ with respect to a $\sigma$-finite measure $\nu(.)$ on $R^d$ and it is known but for the parameter $\theta=(\theta_1,\dots,\theta_k) \in \Theta \in R^k$ which is unknown. The problem is to estimate the parameter $\theta$ based on the observations $\{X^i(t), 0\leq t \leq T\}, i=1,\dots,N.$ 
\vsp
Let $P_\psi$ denote the probability measure generated by the process $X^i$ over the interval $[0,T]$ when $\psi$ is the observe value of the random effect. Let $\psi_0$ be another value of $\psi.$  From the results in Section 2, it follows that, given $\psi,$ the log-likelihood ratio, based on the observation of the process $X^i$ over the interval $[0,T],$ is given by 
$$\ell_T(X^i,\psi)= \int_0^T[Q_{H,\psi}(s)-Q_{H,\psi_0}(s)]dZ_i(s)-\frac{1}{2}\int_0^T[Q^2_{H,\psi}(s)-Q^2_{H,\psi_0}(s)]dw^H_s$$
where the process $Q_{H,\psi}(t), 0\leq t \leq T$ is given by
\be
Q_{H,\psi}(t)= \frac{d}{d<M^H>_t}\int_0^tg_H(s,t)S(X_s^i,\psi) ds, t \geq 0
\ee
as defined by the equation (2.7). Without loss of generality, suppose that the parameter $\psi_0$ satisfies the condition $S(x,\psi_0)=0.$ Then the log-likelihood function, given $\psi,$  can be written in the form
$$\ell_T(X^i,\psi)= \int_0^T Q_{H,\psi}(s)dZ_i(s)-\frac{1}{2}\int_0^TQ^2_{H,\psi}(s)dw^H_s.$$ 
Since the parameter $\phi$ has the probability density $g(\psi,\theta)$ with respect to the $\sigma$-finite measure $\nu(.),$ it follows that the likelihood function based on the observation of the process $X^i$ over the interval $[0,T],$ is given by 
$$\int_{R^d} \exp(\ell_T(X^i,\psi)) \;g(\psi, \theta)d\nu(\psi)= \lambda(X^i,\theta) \;\;\;\mbox{(say)}$$
and the likelihood function based on the processes $\{X^i(t), 0\leq t \leq T\}, i=1,\dots,N$ is given by
$$ L_N(X^1,\dots,X^N;\theta)\equiv \Pi_{i=1}^N\lambda(X^i,\theta).$$
\vsp
Suppose there exists an estimator $\hat \theta_N$ defined by the relation
$$ L_N(X^1,\dots,X^N;\hat \theta_N)= \sup_{\theta \in \Theta}L_N(X^1,\dots,X^N;\theta).$$
The estimator $\hat \theta_N$ is called a maximum likelihood estimator (MLE) of the parameter $\theta.$
Existence and uniqueness of the MLE can be ensured under some conditions on the function $\lambda(f,\theta)$ and the parameter space $\Theta.$
(cf. Prakasa Rao [33]). Under standard regularity conditions on the function $\lambda(f,\theta)$, it can be shown that the MLE $\hat \theta_n$ based on i.i.d. observations is consistent and asymptotically normal as $N \raro \ity$ in the sense that 
$$\hat \theta_N \raro \theta \;\;\mbox{in probability}$$
and
$$\sqrt{N}(\hat \theta_N-\theta) \raro N(0, I^{-1}(\theta)) \;\;\mbox{in law}$$
as $N \raro \ity$ where $I(\theta)$ is the Fisher information matrix defined by
$$I(\theta)= ((I_{i,j}(\theta))), I_{i,j}(\theta)=E[\frac{\partial \lambda (f,\theta)}{\partial \theta_i}\frac{\partial \lambda (f,\theta)}{\partial \theta_j}], i,j=1,\dots,k.$$
See Theorems 16.2 and Theorem 16.3, DasGupta [34] among others for details about the asymptotic theory of maximum likelihood estimators based on independent and identically distributed observations. We omit the details.
\vsp
\newsection{Linear multiplier case}
Suppose that the function $S(x,\psi)=\psi S(x)$ where $\psi\in R$ and the function $S(x)$ is known. We assume that the function $S(x)$ satisfies sufficient conditions so that the system (3.1) has a unique solution. Further suppose that
$$\int_0^TQ_{H,\psi}^2(s)dw^H_s ds <\ity \;\;\mbox{a.s.}$$ 
Observe that the processes $X^i \equiv \{X^i(t), 0 \leq t \leq T\}$ are independent and identically distributed. From the computations given in the Sections 2 and 3, it is easy to check that the likelihood function based on these observed data is given by 
$$\Pi_{i=1}^N\lambda(X^i,\theta)$$
where
\be
\lambda(X^i,\theta)= \int_{R^d}g(\psi,\theta)\exp(\psi \; U_i-\frac{1}{2} \psi^2 \;V_i)\;\;d\nu(\psi),
\ee
\be
U_i= \int_0^T(\frac{d}{dw^H_t}\int_0^t g_H(s,t)S(X^i_s)ds)\;\;dZ_i(s)
\ee
and
\be
V_i= \int_0^T(\frac{d}{dw^H_t}\int_0^tg_H(s,t)S(X^i_s)ds)^2\;\;dw^H_s
\ee
for $i=1,\dots,N.$ Maximizing the likelihood function, we obtain the maximum likelihood estimator of the parameter $\theta.$
\vsp
\NI{\bf Special case:} Suppose that the functions $g(\psi,\theta)$ and the $\sigma$-finite measure $\nu(\psi)$ are such that $g(\psi,\theta)d\nu(\psi)$ is the Gaussian probability density with mean $\mu$ and variance $\sigma_0^2.$ Here $\theta=(\mu,\sigma_0^2).$ Following the computations similar to those given in Dai et al. (2020), it can be shown that, if the parameter $\sigma_0^2$ is known, then the MLE $\hat \mu_N$ of $\mu$ is given by
\be
\hat \mu_N = \sum_{i=1}^N\frac{U_i}{1+\sigma_0^2 V_i}/\sum_{i=1}^N \frac{V_i}{1+\sigma_0^2V_i}.
\ee
If 
$$ 
E(\frac{U_1}{1+\sigma_0^2V_1}) < \ity \;\mbox{and}\;\;E(\frac{V_1}{1+\sigma_0^2V_1})<\ity,
$$
then, by the Strong law of large numbers (SLLN) for independent and identically distributed random variables (i.i.d.), it follows that
\be
\hat \mu_N \raro E(\frac{U_1}{1+\sigma_0^2V_1})/E(\frac{V_1}{1+\sigma_0^2V_1}) = \gamma_0 \;\;\mbox{ (say) a.s. as} \;\;N \raro \ity.
\ee
Applying the SLLN for i.i.d. random variables, it follows that
\be
N^{-1}\sum_{i=1}^N \frac{V_i}{1+\sigma_0^2V_i}\raro E(\frac{V_1}{1+\sigma_0^2V_1})= \beta_0 \;\;\mbox{(say)} \;\;\mbox{a.s. as} \;\;N \raro \ity
\ee
Applying the standard central limit theorem for independent and identically distributed random variables, it can be checked that 
\be
\sqrt{N}(\hat\mu_N- E(\hat \mu_N)) \raro N(0, Var[\frac{U_1}{1+\sigma_0^2V_1}]/\beta_0^2)
\ee
in distribution as $N\raro \ity.$ 
\vsp
If both the parameters $\mu$ and $\sigma_0^2$ are unknown, then the MLEs $\hat \mu_N$ and $\hat \sigma_{0,N}^2$ of $\mu$ and $\sigma_0^2$ are given by the system of equations
\be
\hat \mu_N = (\sum_{i=1}^N\frac{V_i}{1+\hat \sigma_{0,N}^2V_i})^{-1}(\sum_{i=1}^N\frac{U_i}{1+\hat \sigma_{0,N}^2V_i})
\ee
and
\be
\sum_{i=1}^N(\hat \mu_N-\frac{U_i}{V_i})^2\frac{V_i^2}{(1+\hat \sigma_{0,N}^2V_i)^2}= \sum_{i=1}^N \frac{V_i}{1+\hat \sigma_{0,N}^2 V_i}
\ee
but the explicit computation of the estimators $\hat \mu_N$ and $\hat \sigma_{0,N}^2$ and study of their asymptotic properties as $N$ tends to infinity is cumbersome.
\vsp
\NI{\bf Acknowledgment:} This work was supported under the scheme ``INSA Senior Scientist" by the Indian National Science Academy  at the CR Rao Advanced Institute for Mathematics, Statistics and Computer Science, Hyderabad 500046, India. 
\vsp
\NI{\bf References}
\begin{description}

\item {[1]} Mishura, Y. (2008){\it Stochastic Calculus for Fractional Brownian Motion and Related Processes}. Berlin : Springer.

\item {[2]} Prakasa Rao, B.L.S. (2010) {\it Statistical Inference for Fractional Diffusion Processes}, Wiley: Chichester.

\item {[3]} Cheridito, P. (2001) Mixed fractional Brownian motion. {\it Bernoulli}, 7:913-934.

\item {[4]} Mishura, Y. and Shevchenko, G. (2012) Existence and uniqueness of the solution of stochastic differential equation involving Wiener process and fractional Brownian motion with Hurst index $H> 1/2.$ {\it Comput. Math. Appl.} 64: 3217-3227.

\item {[5]} Shevchenko, G. (2014) Mixed stochastic delay differential equations, {\it Theory Probab. Math. Statist.}, 89: 181-195.

\item {[6]} Guerra, J. and Nualart, D. (2008) Stochastic differential equations driven by fractional Brownian motion and and standard Brownian motion, {\it Stochastic Anal. Appl.}, 26: 1053-1075.

\item {[7]} Luis da Silva, Jose., Erraoui, Mohamed and Essaky, El Hassan. (2018) Mixed stochastic differential equations: existence and uniqueness result, {\it J. Theor. Probab.}, 31: 1119-1141.

\item {[8]} Prakasa Rao, B.L.S. (2018) Parameter estimation for linear stochastic differential equations driven by  mixed fractional Brownian motion. {\it Stochastic Anal. Appl.}, 36:767-781.

\item {[9]} Prakasa Rao, B.L.S. (2017) Instrumental variable estimation for a linear stochastic differential equation driven by a mixed fractional Brownian motion. {\it Stochastic Anal. Appl.}, 35: 943-953.

\item {[10]} Prakasa Rao, B.L.S. (2015) Option pricing for processes driven by mixed fractional Brownian motion with superimposed jumps. {\it Probability in the Engineering and Information Sciences}, 29:  589-596.

\item {[11]} Prakasa Rao, B.L.S. (2015) Pricing geometric Asian power options under mixed fractional Brownian motion environment. {\it Physica A}, 446: 92-99.

\item {[12]} Antic, J., Laffont, C.M., and Chafai, D. (2009) Comparison of nonparametric methods in nonlinear mixed effects models, {\it Comput. Statist. Data Anal.}, 53: 642-656.

\item {[13]} Delattre, M., Genon-catalot, V. and Sampson, A. (2012) Maximum likelihood estimation for stochastic  differential equations with random effects, {\it Scand. J. Stat.}, 40: 322-343.

\item {[14]} Ditlevsen, S. and De Gaetano, A. (2005) Mixed effects in stochastic differential equation models, {\it REVSTAT-Statist J.} , 3: 137-153.

\item {[15]} Nie, L. and Yang, M. (2005) Strong consistency of the MLE in nonlinear mixed-effects models with large cluster size, {\it Sankhya}, 67: 736-763.

\item {[16]} Nie, L. (2006) Strong consistency of the maximum likelihood estimator in generalized linear and non-linear mixed-effects models, {\it Metrika}, 63: 123-143.

\item {[17]} Nie, L. (2007) Convergence rate of the MLE in generalized linear and nonlinear mixed-effects models: theory and applications, {\it J. Statist. Plann. Inference}, 137: 1787-1804. 

\item {[18]} Picchini, U., De Gaetano, A., and Ditlevsen, S. (2010) Stochastic differential mixed-effects models, {\it Scand. J. Statist.}, 37: 67-90.

\item {[19]} Picchini, U. and Ditlevsen, S. (2011) Practical estimation of a high dimensional stochastic differential mixed-effects models, {\it Comput. Statist. data Anal.}, 55: 1426-1444.

\item {[20]} El Omari, M., El Maroufy, H. and Fuchs, C. (2019) Nonparametric estimation for fractional diffusion processes with random effects, {\it Statistics}, 53: 753-769.

\item  {[21]} Prakasa Rao, B.L.S. (2021) Nonparametric estimation for stochastic differential equations driven by mixed fractional Brownian motion with random effects, {\it Sankhya Series A} (to appear).

 \item {[22]} Marushkevych, Dmytro. (2016) Large deviations for drift parameter estimator of mixed fractional Ornstein-Uhlenbeck process. {\it Mod.  Stoch. Theory Appl.}, 3:107-117.

\item {[23]} Rudomino-Dusyatska, N. (2003) Properties of maximum likelihood estimates in diffusion and fractional Brownian models. {\it Theor. Probab. Math. Statist.}, 68: 139-146.

\item {[24]} Song, N. and Liu, Z. (2014) Parameter estimation for stochastic differential equations driven by mixed fractional Brownian motion. {\it  Abst. Appl. Anal.}, 2014 Article ID 942307, 6 pp.

\item {[25]} Mishra, M.N. and Prakasa Rao, B.L.S. (2017) Large deviation probabilities for maximum likelihood estimator and Bayes estimator of a parameter for mixed fractional Ornstein-Uhlenbeck type process. {\it Bull. Inform. and Cyber.} 49: 67-80.

\item {[26]} Prakasa Rao, B.L.S. (2009) Estimation for stochastic differential equations driven by mixed fractional Brownian motion. {\it Calcutta Stat. Assoc. Bull.}, 61: 143-153.

\item  {[27]} Miao, Y. (2010) Minimum $L_1$-norm estimation for mixed Ornstein-Uhlenbeck type process, {\it Acta. Vietnam}. 35: 379-386.
 
\item {[28]} El Omari, M., El Maroufy, H. and Fuchs, C. (2019) Statistical inference for fractional diffusion process with random effects  at discrete observations, arXiv: 1912.01463v1 [math.ST] 3 Dec 2019.

\item  {[29]} Cai, C., Chigansky, P. and Kleptsyna, M. (2016) Mixed Gaussian processes. {\it Ann. Probab.}, 44:3032-3075.

\item {[30]} Chigansky, P. and Kleptsyna, M. (2019) Statistical analysis of the mixed fractional Ornstein-Uhlenbeck process, {\it Theory Probab. Appl.}, 63: 408-425.

\item {[31]} Prakasa Rao, B.L.S. (2021) Parametric inference for stochastic differential equations driven by a mixed fractional Brownian motion with random effects based on discrete observations, {\it Stochastic Anal. Appl.}, https:// doi.org/10.1080/07362994.2021.1902352 

\item {[32]} Dai, M., Duan, J., Liao, J., and Wang, X. (2020) Maximum likelihood estimation of stochastic differential equations with random effects driven by fractional Brownian motion, arXiv:2001.01412v! [math.ST] 6 Jan 2020.

\item {[33]} Prakasa Rao, B.L.S. (1987) {\it Asymptotic Theory of Statistical Inference}, Wiley, New York.

\item {[34]} DasGupta, A. (2008) {\it Asymptotic Theory of Statistics and Probability}, Springer, New York.

\end{description}
\end{document}